\documentclass{amsart}

\usepackage{amsmath,amssymb,amscd,verbatim,graphicx}
\usepackage[all]{xy}
\SelectTips{cm}{}
\entrymodifiers={+!!<0pt,\fontdimen22\textfont2>}

\def\vee{*}
\def\too{\longrightarrow}
\def\({\left(}
\def\){\right)}
\def\[{\psi\(\tfrac12B(}
\def\]{)\)}
\def\doplus{\displaystyle\bigoplus}
\def\inv{{-1}}
\def\ad#1{\hat{#1}}

\def\Edge{\mathbb E}
\DeclareMathOperator\image{im}
\def\Sym{\mathrm{Sym}^2}
\def\id{\mathrm{id}}

\def\CCC{\mathbb C}
\def\RRR{\mathbb R}

\def\DD{D}

\def\CC{C}

\def\GL{{\mathit{GL}}}
\def\gl{{\mathfrak{gl}}}

\def\t#1#2#3{\tau(l_{#1},l_{#2},\ldots,l_{#3})}
\def\q#1#2#3{q_{{#1},{#2},\ldots,{#3}}}
\def\bq#1#2#3{\tilde q_{{#1},{#2},\ldots,{#3}}}
\def\T#1#2#3{T_{{#1},{#2},\ldots,{#3}}}
\def\Es#1#2#3{S_{{#1},{#2},\ldots,{#3}}}

\def\bT#1#2#3{\tilde T_{{#1},{#2},\ldots,{#3}}}
\def\K#1#2#3{K_{{#1},{#2},\ldots,{#3}}}
\def\bK#1#2#3{\tilde K_{{#1},{#2},\ldots,{#3}}}

\def\V{V}

\def\tK{\tau^{\mathrm{Kash}}}
\def\TK{T^{\mathrm{Kash}}_{1,2,3}}

\def\qK{q^{\mathrm{Kash}}_{1,2,3}}

\def\<{\left<}
\def\>{\right>}

\def\WF{W(F)}

\theoremstyle{plain}
\newtheorem{thm}{Theorem}[section]
\newtheorem{prop}{Proposition}
\newtheorem{lemma}[prop]{Lemma}

\newtheorem*{cor}{Corollary}

\theoremstyle{definition}

\theoremstyle{remark}
\newtheorem*{rem}{Remark}
\newtheorem{remnum}{Remark}
\newtheorem*{Case1}{Transverse Case}
\newtheorem*{Case2}{General Case}

\def\ZnZ{\mathbb Z/n\mathbb Z}

\title{The Maslov index as a quadratic space}
\author{Teruji Thomas}
\address{Merton College, Oxford University, Oxford OX1 4JD, UK}%
\email{Joaquin.Thomas@aya.yale.edu}%
\thanks{Partially supported by the University of Chicago's  VIGRE grant,  DMS-9977134, and by NSF grant DMS-0401164. \\
{\bf Note:} A version of this paper is published in {\it Math. Res. Lett.} {\bf 13}, vol 6. (2006). In the published version, sections 9--11 are omitted for reasons of space, and section 3 is slightly rewritten.}

\begin{document}

\begin{abstract}
Kashiwara defined the Maslov index (associated to a collection of Lagrangian subspaces of a symplectic vector space over a field $F$) as a class in the Witt group $W(F)$ of quadratic forms. We construct a canonical quadratic vector space in this class and show how to understand the basic properties of the Maslov index without passing to $W(F)$---that is, more or less, how to upgrade Kashiwara's equalities in $W(F)$ to canonical isomorphisms between quadratic spaces. We also show how our canonical quadratic form occurs naturally in the context of the Weil representation. The quadratic space is defined using elementary linear algebra. On the other hand, it has a nice interpretation in terms of sheaf cohomology, due to A. Beilinson.
\end{abstract}

\maketitle

\section{Introduction.} Let $F$ be a field of characteristic not 2. Let $\WF$ denote the Witt group of quadratic spaces over $F$---we will use the term {\it quadratic space} to mean a finite dimensional $F$-vector space with a non-degenerate symmetric bilinear form.

\subsection{}\label{Kashconstruction}
Suppose given a  vector space $\V$ over $F$ with a symplectic form $B$ and some Lagrangian subspaces $l_1,\ldots,l_n$ of $\V$, indexed by $\ZnZ$. 

To this data Kashiwara associated a class $\tau(l_1,\ldots,l_n)\in\WF$ called {\it the Maslov index} of the Lagrangians (see the appendices in \cite{LV} and \cite{KS}, and section \ref{sec:Kash} of this article). 
In this article we answer the following questions:
\begin{itemize}
\item[(i)] How can one represent $\tau(l_1,\ldots,l_n)$ as the class of a canonically defined quadratic space?
\item[(ii)] How can one upgrade the basic equalities satisfied by $\tau$ to canonical isomorphisms between quadratic spaces?
\end{itemize}
By `basic equalities' we mean dihedral symmetry, i.e. 
\begin{equation}\label{Kash:dihedral}
\tau(l_1,l_2,\ldots,l_n)=\tau(l_2,l_3,\ldots,l_{n},l_1)=-\tau(l_n,l_{n-1},\ldots,l_1)
\end{equation}
and the chain condition, i.e.
\begin{equation}\label{Kash:cocycle}
\tau(l_1,l_2,\ldots,l_n)=\tau(l_1,l_2,\ldots,l_k)+\tau(l_1,l_k,\ldots,l_n)
\end{equation}
for any $k\in\{3,\ldots,n-1\}$. 
We also mean that when $F$ is a local field (e.g. $F=\mathbb R$), 
\begin{equation}\label{Kash:constant}
\substack{\mbox{\it $\tau(l_1,\ldots,l_n)$ is locally constant in $l_1,\ldots,l_n$} \\
\mbox{\it if the dimension of $l_i\cap l_{i+1}$ is fixed for each $i\in\ZnZ$.} 
}
\end{equation}

\subsection{} In section \ref{sec:definitions} we answer question (i), constructing  a quadratic space denoted $\T12n$ (with bilinear form $\q12n$). We re-define the Maslov index $\t12n$ to be the class of $\T12n$ in $W(F)$, and in section \ref{sec:Kash} we verify that our Maslov index is the same as Kashiwara's. 

\label{sec:dodo} We will give a concrete description of the quadratic space in section \ref{concrete}, but abstractly we may say that the vector space $\T12n$ is the cohomology $H^0(C)$ of a certain complex $C$ \eqref{complex2}; the form $\q12n$ arises from the quasi-isomorphism of $C$ and its dual $C^*$. This point of view is explained in section \ref{sec:nondeg}, where we also give a formula for the dual form on $\T12n^*$.

 A sheaf-theoretic construction of $(\T12n,\q12n)$ is described in \ref{intro:Beilinson}. The reader chiefly interested in this interpretation can proceed directly to \ref{intro:Beilinson} and then section \ref{sec:Beilinson}.

\begin{rem} Having written this paper, we noticed that another answer to question (i) is proposed in \cite{CLM}, section 12,  using topological methods when $F=\mathbb R$.  However, one can find a counter-example to their formula for $n=5$, $\dim V=2$. Their method rather leads to our formula \eqref{eq:explicit}.
\end{rem}

\subsection{} In sections \ref{sec:dihedral}, \ref{sec:cocycle}, \ref{sec:constant} we answer question (ii).

 \subsubsection{} The dihedral symmetry (\ref{Kash:dihedral}) will be realized by canonical identifications 
\begin{equation*}\label{dihedral}
(\T12n,\q12n)=(\T{2}3{n,1},\q{2}3{n,1})=(\T{n}{n-1}1,-\q{n}{n-1}1)
\end{equation*}
as described in section \ref{sec:dihedral}.

\subsubsection{} For the chain condition (\ref{Kash:cocycle}), suppose first that $l_1\cap l_k=0$. In section \ref{sec:cocycle} we describe a canonical isometric isomorphism
\begin{equation*}
\begin{CD}
\T12k\oplus\T1kn @>\cong>> \T12n.
\end{CD}
\end{equation*}
To treat the case $l_1\cap l_k\neq0$,  we use the notion of {\it quadratic subquotient}:  if $T$ is a  quadratic space and $I$ an isotropic subspace, then $I^\perp/I$ is again a quadratic space, called `the quadratic subquotient of $T$ by $I$.'  

Without assuming $l_1\cap l_k=0$,  we identify 
 $\T12k\oplus\T1kn$ with a quadratic subquotient of $\T12n$, so the following well known lemma shows that (\ref{Kash:cocycle})
holds without any assumptions. 

\begin{lemma}{\label{prop:criterion} If $S$ is a quadratic subquotient of $T$ then $S$ and $T$ have the same class in $\WF$. }\end{lemma}
\begin{proof}{ Suppose that $S=I^\perp/I$. Choose a linear complement $M$ to $I^\perp$ in $T$. Then $M+I$ (a direct sum) is a hyperbolic summand of $S$, and $(M+I)^\perp\subset I^\perp$ maps isometrically onto $I^\perp/I$.}\end{proof}

\subsubsection{} As for \eqref{Kash:constant}, it is not true in general that the isomorphism class of $\T12n$ is locally constant under the condition described. It would be enough to show that $\dim\T12n$ is locally constant (see Lemma \ref{prop:dimonly}), but in fact (see \eqref{eq:dim}), $\dim\T12n$ depends not only on the dimensions of $l_i\cap l_{i+1}$ but also on $\dim (l_1\cap\cdots\cap l_n)$. 

In section \ref{sec:constant} we represent $\T12n$ as a quadratic subquotient of another quadratic space $\bT12n$ whose dimension depends only on the dimensions $\dim l_i\cap l_{i+1}$. Property \eqref{Kash:constant} then follows from Lemma \ref{prop:criterion} above.

\begin{rem} As already mentioned in \ref{sec:dodo}, there is a quasi-isomorphism $\Phi\colon C\to C^*$. To define $\bT12n$ we factor $\Phi$ as 
$$
\begin{CD}
C @>\alpha>> D @>\tilde\Phi>>  D^* @>\alpha^*>> C^*
\end{CD}
$$
where $\alpha$ is a quasi-isomorphism and $\tilde\Phi$ an isomorphism. The construction of such a factorization (see \ref{sec:factor}) is quite general and is useful in other situations (see e.g. \cite{Wa}, \cite{So}, \cite{Ke}). 
\end{rem}

\subsection{} \label{intro:Beilinson} Here is a sheaf-theoretic interpretation due to A.~Beilinson (private communication). He interprets $\T12n$ as $H^1(D,P)$, where $D$ is the filled $n$-gon and $P$ is the following subsheaf of the constant sheaf $V_D$ with fibre $\V$.  Cyclically label the vertices of $D$ by $\ZnZ$ and give $D$ the corresponding orientation. Let $U$ be the interior of $D$. Then $P$ has fibre $\V$ on $U$, $l_i$ on the edge $(i,i+1)$, and $l_{i-1}\cap l_{i}$ at the vertex $i$. 

Write $j_U\colon U\to D$ for the inclusion. Let $F_U$ be the constant sheaf on $U$ with fibre $F$.  The symplectic form on $\V$ induces a map $P\otimes P\to j_{U,!}F_U$, and thereby
\begin{equation}
\label{cup}\cup\colon\Sym H^1(D,P)\to H^2(D,j_{U,!}F_U)=F.
\end{equation}
 In section \ref{sec:Beilinson}, which can mostly be read independently from the rest of the work, we show that \eqref{cup} is non-degenerate and in fact is {\it minus} our original bilinear form.

\subsubsection{} From this perspective, dihedral symmetry \eqref{Kash:dihedral} is already clear, since the cyclic symmetry is manifest, while reversing the order of the Lagrangians reverses the orientation of $D$ and so changes the identification $H^2(D,j_!F_U)=F$ by a sign.

\subsubsection{} In \ref{sec:BeilinsonChain} we explain how the chain condition {\eqref{Kash:cocycle}} may be proved via a ``bordism'' between three polygons, in analogy to the proof of the additivity of the index of manifolds.

\subsection{} Finally, in sections \ref{sec:gamma}--\ref{sec:final} we show, when $F$ is a finite or local field, how our quadratic form occurs naturally in the context of the Weil representation, which originally motivated this work. The relationship between the Maslov index and the Weil representation is well known (see \cite{LV},\cite{Li},\cite{Pe},\cite{Ra},\cite{Sou},\ldots), with $n=3$ Lagrangians being the key case. Our definitions allow a direct approach for any $n$. 

In section \ref{sec:final} we also describe the self-dual measure on $\T12n$.

\subsection{} I am grateful to V.~Drinfeld for suggesting this subject, and to him and D.~Arin\-kin for much useful advice. I also thank A.~Beilinson for explaining his sheaf-theoretic reformulation, and M. Kamgarpour and B. Wieland for many interesting discussions.

\section{The Quadratic Space.}
\label{sec:definitions}
\subsection{Preliminaries. }\label{sec:prelim} Our Lagragians $l_1,\ldots,l_n$ are indexed by $\ZnZ$.  Think of $\ZnZ$ as the vertices of a graph whose set $\Edge$ of edges consists of pairs of consecutive numbers $\{i,i+1\}$, $i\in\ZnZ$. So the graph looks like an $n$-sided polygon. 

An element $$v=(v_{\{i,i+1\}})\in\doplus_{\{i,i+1\}\in\Edge}\V$$
may be thought of as a function $a\colon \Edge\to\V$.  We can form the `derivative' 
\begin{equation}\label{derivative}
\partial v=(\partial v_i)\in\doplus_{i\in\ZnZ} \V \qquad\qquad \partial v_i=v_{\{i,i+1\}}-v_{\{i-1,i\}}.
\end{equation}
Conversely, given $w=(w_i)\in\bigoplus_{i\in\ZnZ} \V,$ there is the obvious notion of an `anti-derivative' $\ad{w}$
\begin{equation}\label{eq:anti}
\ad{w}=(\ad{w}_{\{i,i+1\}})\in\!\!\!\doplus_{\{i,i+1\}\in\Edge}\V \qquad\mbox{such that}\quad \partial(\ad{w})=w.\end{equation}
 An anti-derivative exists so long as 
$\sum_{i\in\ZnZ} w_i=0$
in which case $\ad{w}$ is unique up to adding a constant function.

\subsection{Definitions.} \label{concrete}

\subsubsection{\bf Definition}
  Let $\K12n$ be the kernel of the natural summation 
\begin{equation}\label{eq:defK}
\K12n=\ker  \left[\bigoplus_{i\in\ZnZ} l_i\overset\Sigma\longrightarrow\V\right].
\end{equation} 

Any $w=(w_i)\in\K12n$ has an anti-derivative $\ad{w}\in\bigoplus_{\{i,i+1\}\in\Edge}\V$, as in \ref{sec:prelim}.

\subsubsection{\bf Definition} { Define a bilinear form $\q12n$ on $\K12n$ by the formula
\begin{equation}\label{eq:defq}
\q12n(v,w)=\sum_{i\in\ZnZ} B(v_i, \ad w_{\{i,i+1\}})
\end{equation}
for any choice of anti-derivative $\ad{w}$ (it is simple to check that the right-hand side of \eqref{eq:defq} is independent of this choice.)}

\begin{remnum}
{The bilinear form may equivalently be defined by
\begin{equation}\label{alternate}
\q12n(v,w)=\sum_{i\in\ZnZ} B(v_i, \ad w_{\{i-1,i\}}).
\end{equation}
Indeed, the difference  is
$$\sum_{i\in\ZnZ} B(v_i,\ad w_{\{i,i+1\}}-\ad w_{\{i-1,i\}})=\sum_{i\in\ZnZ} B(v_i,w_i),$$
and $B(v_i,w_i)=0$ since $v_i,w_i$ lie in the same Lagrangian $l_i$.}
\end{remnum}

\begin{remnum} In the definition \eqref{eq:defq} of $\q12n$, one may concretely choose
$\ad w_{\{i,i+1\}}=\sum_{j=1}^i w_j$,
in which case the definition takes the simple form
\begin{equation}\label{eq:explicit}
\q12n(v,w)=\sum_{i\geq j\geq 1}^n B(v_i,w_j)=\sum_{i>j>1}^n B(v_i,w_j).
\end{equation}
We often use this version for calculations, but the natural symmetries of $\q12n$ are obscured.
\end{remnum}

\begin{prop}{The bilinear form $\q12n$ is symmetric.}\end{prop}
\begin{proof}{``Summation by parts.'' Explicitly, 
\begin{equation*}
\begin{aligned}
\q12n(v,w)&=\sum_{i\in\ZnZ} B(v_i,\ad w_{\{i,i+1\}})\\
&=\sum_{i\in\ZnZ} B(\ad v_{\{i,i+1\}}-\ad v_{\{i-1,i\}},\ad w_{\{i,i+1\}})\\
&=\sum_{i\in\ZnZ} B(\ad w_{\{i,i+1\}}-\ad w_{\{i-1,i\}},\ad v_{\{i-1,i\}})\\
&=\sum_{i\in\ZnZ} B(w_i,\ad v_{\{i-1,i\}})=\q12n(w,v).
\end{aligned}
\end{equation*}
The last equality is \eqref{alternate}.}\end{proof}

\subsubsection{\bf Definition of the Quadratic Space $(\T12n,\q12n)$}\label{TheDef}{ 
\begin{equation}
\T12n=\K12n\big/\ker\q12n.
\end{equation} 
The induced non-degenerate bilinear form on $\T12n$ will still be called $\q12n$.}

\smallskip
Let us give an explicit description of $\T12n$. The derivative \eqref{derivative} restricts to give a map
\begin{equation}\label{deriv2}
\partial\colon\doplus_{\{i,i+1\}\in\Edge}l_i\cap l_{i+1}\to \K12n.
\end{equation}
It is clear from the definition  \eqref{eq:defq} that  $\ker\q12n\supset\image\partial$.

\begin{prop} 
\label{prop:nondeg}We have $\ker\q12n=\image\partial.$
In other words, $\T12n$ is the cohomology $H^0(C)$ at the center term of the complex
\begin{equation}\label{complex2}
\begin{CD}
C=[\doplus_{\{i,i+1\}\in\Edge} l_i\cap{ l_{i+1}} @>\partial>> \doplus_{i\in\ZnZ}  l_i @>\Sigma>>  \V]
\end{CD}
\end{equation}
which we consider to lie in degrees $-1,0,1.$
\end{prop}

The proof will be given in section \ref{sec:nondeg}, where we also show that $\q12n$ is induced by a quasi-isomorphism between $C$ and $C^*$, and give a formula for the dual form on $\T12n^*$.

\begin{cor}
\begin{equation}\label{eq:dim}
\dim \T12n=(n-2)\frac{\dim\V}2-\sum_{i\in\ZnZ}\dim l_i\cap l_{i+1}+2\dim\bigcap_{i\in\ZnZ} l_i.
\end{equation}
\end{cor}

\begin{rem} The discriminant of $\q12n$ has also be determined; see \cite{Th}.
\end{rem}

\subsubsection{\bf Definition}{The symbol $\t12n$ denotes the class of the quadratic space $(\T12n,\q12n)$ in $\WF$, called {\it the Maslov index} of $l_1,\ldots,l_n$.}

\bigskip
In section \ref{sec:Kash} we will verify that this Maslov index equals Kashiwara's.

\section{ Proof of Proposition \ref{prop:nondeg}; \\ The Dual Form; Homotopy Equivalence of $C$ and $C^*$. }
\label{sec:nondeg}

In this section we give an algebraic proof of Proposition \ref{prop:nondeg} (a sheaf-theoretic argument is given in section \ref{sec:Beilinson}). The proof allows us to write down a formula for the dual form on $\T12n^*$ in \ref{sec:dualform}.
  It also implies that the complexes $C$ and $C^*$ are isomorphic in the derived category of complexes of vector spaces, this isomorphism inducing $\q12n$ on $H^0(C)$.  In \ref{sec:htp} we give an explicit, though non-canonical, symmetric quasi-isomorphism $C\to C^*$.

\subsection{Proof of Proposition \ref{prop:nondeg}.}\label{nondegproof} Let $\Phi$ be the linear map $\Phi\colon H^0(C)\to H^0(C^*)=H^0(C)^*$ such that $$\q12n(v,w)=\<v,\Phi w\>.$$  We want to show that $\Phi$ 
is an isomorphism.

The situation is expressed by the following commutative diagram,
in which the top row is $C$ and the bottom row  $C^*$; every row is a complex.
\begin{equation}\label{seq:bd1}
\xymatrix{
& {\doplus_{\{i,i+1\}\in\Edge} 
l_i\cap l_{i+1}} 
 \ar[r]^(0.53){\partial}="x"\ar[d]^-{\mathrm{incl}}
  & {\doplus_{i\in\ZnZ} l_i}  \ar[r]^(0.67){\Sigma}\ar[d]^-{\mathrm{incl}}
  & {\V}\ar@{=}[d]
  \\
  & {\doplus_{i\in\ZnZ} l_i} \ar[r]^(0.49){\partial}="y"\ar[d]^{\mathrm{incl}}
  & {\doplus_{\{i,i+1\}\in\Edge} l_i+l_{i+1}}\ar[r]^(0.66){\Sigma}\ar[d]^{\mathrm{incl}}
  & {\V}\ar@{=}[d]\\
 {\V} \ar[r]^(0.4){\mathrm{diag}} \ar[d]^{\mathrm{proj}}
  & {\doplus_{i\in\ZnZ} \V} \ar[r]^!{"x";"y"}{\partial}\ar[d]^{\mathrm{proj}}
  & {\doplus_{\{i,i+1\} \in\Edge} \V}\ar[r]^(0.67){\Sigma}\ar[d]^{\mathrm{proj}}
  & {\V} \\
  {\V}\ar[d]^{\beta} \ar[r]^(0.4){\mathrm{diag}}
  & {\doplus_{i\in\ZnZ} V/l_i}\ar[d]^{\beta}\ar[r]^!{"x";"y"}{\partial}
  & {\doplus_{\{i,i+1\}\in\Edge}V/(l_i+l_{i+1})} \ar[d]^{-\beta}\\
  {\V^*} \ar[r]^(0.4){\mathrm{diag}}
  & {\doplus_{i\in\ZnZ} l_i^*}\ar[r]^!{"x";"y"}{\partial^*}
  & {\doplus_{\{i,i+1\}\in\Edge}(l_i\cap l_{i+1})^*}.
}
\end{equation}
Here ``$\mathrm{incl}$'' and ``$\mathrm{proj}$'' mean summand-by-summand inclusion and projection, and $\beta\colon x\mapsto B(-,x)$. 
By definition \eqref{eq:defq} of $\q12n$, $\Phi\colon H^0(C)\to H^0(C^*)$ factors as 
\begin{equation}\label{fPhi}\Phi=\beta\circ\mathrm{proj}\circ\partial^\inv\circ\mathrm{incl}\circ\mathrm{incl}.\end{equation}

The map ``$\mathrm{incl}$'' between the first two rows of \eqref{seq:bd1} is a quasi-isomorphism of complexes: it is injective with acyclic cokernel.
  Denote by $W$ the cohomology at the center term of the second row;
we obtain an isomorphism
$\mathrm{incl}\colon H^0(C) \to W$.

 The second, third, and fourth rows of \eqref{seq:bd1} form a short exact sequence of complexes. Let $W'$ denote the cohomology at the center term of the fourth row. Since the third row is exact, the 
 boundary map $\delta=\mathrm{incl}^\inv\circ\partial\circ\mathrm{proj}^\inv\colon W'\to W$ is an isomorphism.  

Finally, the map $\beta$ between the last two rows of \eqref{seq:bd1} is an isomorphism. 
Thus $\Phi$ factors through isomorphisms
$$\Phi\colon H^0(C) \overset{\mathrm{incl}}\too  W \overset{\delta^\inv}\too W' \overset{\beta}\too H^0(C^*).
$$\qed

\subsection{The Dual Form}\label{sec:dualform} Now let us record a formula for the `dual' quadratic form $\q12n^*$ on $\T12n^*$, defined, in the notation of \ref{nondegproof}, by $\q12n^*(x,y)=\<x,\Phi^\inv y\>$.  

\subsubsection{}\label{defEs}  Define $$\Es12n:=\{x\in\bigoplus_{i\in\ZnZ}V/l_i\,|\,x_{i+1}-x_{i}\in l_i+l_{i+1}\}.$$ $\T12n^*$ is a quotient of $\Es12n$: the map $\Es12n\to\T12n^*$ is given by the last two rows of \eqref{seq:bd1}, in which  $\Es12n=\ker \partial$ (in the fourth row) and $\T12n^*=\ker(\partial^*)/\image(\mathrm{diag})$. We will, more precisely, describe
 $\q12n^*$ pulled back to $\Es12n$.

\def\rii{\varepsilon_{i,i+1}}

\subsubsection{}\label{sec:r} 
 Suppose given $x_i,x_{i+1}\in V$ with $x_{i+1}-x_{i}\in l_i+l_{i+1}$. Define a functional $\rii(x_i,x_{i+1})$ on $l_i+l_{i+1}$ as follows. For $v\in l_i+l_{i+1}$, write $v=a+b$, with $a\in l_i,$ $b\in l_{i+1}$. Then $$\<\rii(x_i,x_{i+1}),v\>:=B(a,x_i)+B(b,x_{i+1}).$$
It is easy to verify that this quantity is independent of the choice of $a,b$.
 
\begin{prop}\label{lem:dualform} As a bilinear form on $\Es12n$,
\begin{equation}\label{eq:dualform} 
\q12n^*(x,y)=\sum_{i\in\ZnZ} \<\rii(x_i,x_{i+1}),y_{i+1}-y_{i}\>.
\end{equation}
\end{prop}

\subsection{An Explicit Quasi-Isomorphism} \label{sec:htp} Here is one particular quasi-isomorph\-ism $\Phi\colon C\to C^\vee$ inducing $\q12n$ on $H^0(C)$:
\begin{equation}\label{seq:quis}
\xymatrix{
   {\doplus_{\{i,i+1\}\in\Edge} l_i\cap{ l_{i+1}}} 
               \ar[r]^-{\partial}\ar[d]^(0.55){\Phi_\inv}
 & {\doplus_{i\in\ZnZ}  l_i} 
               \ar[r]^-{\Sigma}\ar[d]^(0.55){\Phi_0} 
 & \V
               \ar[d]^(0.55){\Phi_1}\\
   {\V^*} \ar[r]^(0.4){\Sigma^*}
 & {\doplus_{i\in\ZnZ}  l_i^*}\ar[r]^(0.4){\partial^*}
 & {\doplus_{\{i,i+1\}\in\Edge} (l_i\cap{ l_{i+1}})^*}
}
\end{equation}
where
\begin{equation}\label{explicitquis}
\begin{aligned}
\<\Phi_{-1}(a),v\>&=\<a,\Phi_{1}(v)\>=B(a_{\{n,1\}},v) \\
\<\Phi_{0}(a),b\>&=\<a,\Phi_0(b)\>=\tfrac12\sum_{i\geq j\geq1}^n\(B\(a_i,b_j\)+B\(b_i,a_j\)\).
\end{aligned}
\end{equation}
The fact that $\Phi$ induces $\q12n$ follows from formula \eqref{eq:explicit}.

\begin{rem} $\Phi=\Phi^*$. The isomorphism $$H^\inv(C)\oplus H^0(C)\oplus H^1(C)\longrightarrow H^\inv(C)^*\oplus H^0(C)^*\oplus H^1(C)^*$$ implicit in diagram \eqref{seq:bd1} is symmetric, so any $\Phi\colon C\to C^*$ inducing it may be symmetrized by $\Phi\mapsto \tfrac12(\Phi+\Phi^*)$.
\end{rem}

\section{Dihedral Symmetry (1).}\label{sec:dihedral}

The space $\K12n$ can be canonically identified with $\K 23{n,1}$ and $\K n{n-1}1$, as is obvious from definition \eqref{eq:defK}. 

\begin{prop}{Under these identifications, $q_{2,3,\ldots,n,1}=\q12n=-q_{n,n-1,\ldots,1}.$ }\end{prop}

\begin{proof}{The first equality (cyclic symmetry) is obvious from the definition \eqref{eq:defq} of $\q12n$. Reversing the order of the Lagrangians is equivalent to replacing $\partial$ by $-\partial$, and therefore $\ad{w}$ by $-\ad{w}$ in equation \eqref{eq:defq}.}\end{proof}

\section{The Chain Condition (2).}\label{sec:cocycle}

\begin{prop}\label{prop:cocycle} Fix $k\in\{2,\ldots,n\}$. \begin{itemize} 
\item[1.] If $l_1\cap l_k=0$ then $\T12k\oplus\T1kn\cong\T12n$ isometrically.
\item[2.] Without conditions, $\T12k\oplus\T1kn$ is a quadratic subquotient of $\T12n.$ 
\end{itemize}
Therefore, by Lemma \ref{prop:criterion}, $\t12k + \t1kn=\t12n$.
\end{prop}

\begin{proof}{The equivalence is induced by the natural map 
\begin{equation*}
\begin{CD}
\doplus_{i\in\{1,2,\ldots,k\}}\!\!\!\!\!l_i\,\oplus\!\!\!\!\doplus_{i\in\{1,k,\ldots,n\}}\!\!\!\!\!l_i @>s>> \doplus_{i\in\{1,2,\ldots,n\}}\!\!\!\!\! l_i
\end{CD}
\end{equation*}
that is the identity on each summand.
More precisely, consider the map of short exact sequences 
\begin{equation*}
\xymatrix@R=0.5in{
{l_1\oplus l_k}
     \ar[r]\ar[d]^(0.55){\Sigma} 
 &  {\!\!\!\!\doplus_{i\in\{1,2,\ldots,k\}}\!\!\!\!\!l_i\,\oplus\!\!\!\!\doplus_{i\in\{1,k,\ldots,n\}}\!\!\!\!\!l_i}
     \ar[r]^(0.6){s}\ar[d]^(0.55){\Sigma\oplus\Sigma}
 &  {\doplus_{i\in\{1,2,\ldots,n\}}\!\!\!\!\! l_i} 
     \ar[d]^(0.55){\Sigma}\\
{\V} \ar[r]^-{(\id,-\id)}
 &  {\V\oplus\V}
     \ar[r]^(0.6){s} 
 &  {\V}
}
\end{equation*}
The snake lemma gives an exact sequence
\begin{equation}\label{seq:longexact}
\begin{CD}
 l_1\cap l_k @>(r,-r)>> \K12k\oplus\K1kn@>s>> \K12n @>\delta>> \V/(l_1+l_k).
\end{CD}\end{equation}
Here $r=(-\id,\id)\colon l_1\cap l_k\to l_1\oplus l_k$ and the boundary $\delta$ is given by
\begin{equation}{\label{eq:imagec}
\delta(v)=\sum_{i=1}^k v_i\mod(l_1+l_k).
}\end{equation}

\begin{lemma}{The map $s\colon \K12k\oplus\K1kn\to \K12n$ is an
isometry.}\end{lemma}

\begin{proof}{ One sees immediately from the explicit formula (\ref{eq:explicit}) that $s$ restricted to each summand is an isometry, and that for any $(v,w)\in \K12k\oplus \K1kn$ we have $\q12n(s(v,0),s(0,w))=0$. No more is required.
}\end{proof}

\begin{Case1}
{If $l_1\cap l_k=0$ then, by the exactness of (\ref{seq:longexact}), $s$ is an isometric isomorphism. Passing to the non-degenerate quotients  establishes part 1 of the proposition.}
\end{Case1}

\begin{Case2} 
Consider the map $r\colon l_1\cap l_k\to\K12k$ as in (\ref{seq:longexact}). By Proposition \ref{prop:nondeg}, the image lies in $\ker\q12k$, 
so, since $s$ is an isometry,  $\image s\circ r$ is isotropic in $\K12n$ and $\image s\subset(\image s\circ r)^\perp$. 

\begin{lemma}{ The image of $s$ is exactly $(\image s\circ r)^\perp.$}\end{lemma}
\begin{proof}{
According to formula (\ref{eq:explicit}), if  $v\in l_1\cap l_k$ and $w\in \K12n$, then
$$
\q12n(s\circ r(v),w) =\sum_{j=1}^k B(v,w_j).$$
This quantity vanishes for all $v\in l_1\cap l_k$ if and only if $\sum_{i=1}^k w_k$ lies in $l_1+l_k$; according to (\ref{eq:imagec}) this just means $\delta(w)=0$ or equivalently $w\in\image s$. 
}\end{proof}

Let $I$ denote the image of $\image s\circ r$ in $\T12n$. We have constructed a surjective isometry
$$s\colon \K12k\oplus\K1kn\longrightarrow I^\perp/I$$
and therefore
$$s\colon \T12k\oplus\T1kn\overset\cong\longrightarrow I^\perp/I.$$
\end{Case2}}
\end{proof}

\section{Local Constancy (3).}\label{sec:constant}

In this section, the ground field $F$ is a local field. 
\subsection{}
 In order to prove \eqref{Kash:constant} we will use the following well known fact:
\begin{lemma} \label{prop:dimonly} Let $Q(X)$ be the space of non-degenerate symmetric bilinear forms on a fixed vector space $X$. Then the natural $\GL(X)$ action on $Q(X)$ has open orbits.
In other words, in a continuous family of quadratic spaces of fixed rank, the isomorphism class is locally constant. \end{lemma}
\begin{proof}
It suffices to check that for any $q\in Q(X)$ the map $\lambda\colon \GL(X)\to Q(X)$ given by $\lambda(g)(x,x)=q(gx,gx)$ has surjective differential at $1\in\GL(X)$. The differential $\lambda_*\colon \gl(X)\to \Sym(X^*)$ is given by $\lambda_*(A)(x,x)=2q(Ax,x)$; it is surjective since $q$ is non-degenerate.
\end{proof}

  We define in this section another quadratic space $\bT12n$ whose class in $\WF$ is the same as that of $\T12n$, but whose dimension is locally constant in $l_1,\ldots,l_n$ 
when the dimensions $\dim l_i\cap l_{i+1}$ are fixed. Then property \eqref{Kash:constant} follows from Lemma \ref{prop:dimonly}. Note that, by \eqref{eq:dim}, the dimension of $\T12n$ itself varies with $\dim (l_1\cap\cdots\cap l_n)$.

\begin{rem}
One could use $\bT12n$ to define the Maslov index, but then the dihedral symmetry \eqref{Kash:dihedral} would not be immediately clear.
\end{rem}

We will first give a concrete formula \eqref{eq:explicitb} for the quadratic form $\bq12n$ on $\bT12n$, but in section \ref{sec:factor} we describe its origin in a general construction.

\subsection{\bf Definition}
 Define a symmetric bilinear form $\bq12n$ on 
\begin{equation}\bK12n=\V^*\oplus \doplus_{i\in\ZnZ}l_i
\end{equation}
by 
\begin{equation}\label{eq:explicitb}
\bq12n(v\oplus a,v\oplus a)=2\sum_{i\in\ZnZ}\<v, a_i\>+\sum_{i\geq j\geq 1}B(a_i,a_j).
\end{equation}
Define 
\begin{equation}\bT12n=\bK12n/\ker\bq12n
\end{equation}
 and again write $\bq12n$ for the non-degenerate form on $\bT12n.$

To describe $\bT12n$ explicitly, let 
$$\tilde\partial=(\Phi_{-1},\partial)\colon\doplus_{\{i,i+1\}\in\Edge}l_i\cap l_{i+1}\to\bK12n.$$
Here $\Phi_{-1}\colon\bigoplus(l_i\cap l_{i+1})\to \V^*$ as in (\ref{seq:quis},\ref{explicitquis}) and $\partial\colon\bigoplus(l_i\cap l_{i+1})\to \bigoplus l_i$ as in \eqref{deriv2}.

\begin{prop} \label{prop:cc} With notation as above:
\nobreak\begin{itemize}
\item[1.] We have $\ker\bq12n=\image\tilde\partial$, and $\tilde\partial$ is injective. Therefore 
$$\dim \bT12n=(n+2)\frac{\dim\V}2-\sum_{\{i,i+1\}\in\Edge}\dim l_i\cap l_{i+1}.$$

\item[2.]
$\T12n$ is a quadratic subquotient of $\bT12n$ by the image of $\V^*\subset\bK12n$. 
Therefore, by Lemma \ref{prop:criterion}, $\T12n$ and $\bT12n$ have the same class in $\WF$. 
\end{itemize}
\end{prop}

\begin{proof}[Proof of Proposition \ref{prop:cc}] 
It is simple to check Proposition $\ref{prop:cc}$ directly, 
making reference to the quasi-isomorphism (\ref{seq:quis},\ref{explicitquis}); in the notation there,
$$\bq12n(v\oplus a,w\oplus b)=\<v,\Sigma b\>+\<a,\Phi_0b\>+\<a,\Sigma^*w\>.$$
 Otherwise, one can feed that quasi-isomorphism into the following general construction.

\subsection{}{\label{sec:factor} {\it Let  $\CC$ be a complex of finite-dimensional vector spaces. Suppose given a symmetric quasi-isomorphism $\Phi\colon \CC\to\CC^\vee$. Then we construct another such complex $\DD$ and a commutative diagram
\begin{equation}\label{eq:factor}
\begin{CD}
\CC @>\alpha>> \DD \\
@V\Phi VV @V\tilde\Phi VV \\
\CC^* @<\alpha^*<< \DD^* \\
\end{CD}
\end{equation}
where $\alpha$ is a quasi-isomorphism and $\tilde\Phi$ a symmetric isomorphism.}}

\begin{rem} Instead of vector spaces, one can consider projective modules over an associative ring with anti-involution. The construction still works, if one everywhere replaces ``quasi-isomorphism'' with ``homotopy equivalence'' and ``acyclic'' with ``homotopy equivalent to zero.''

In fact, \cite[Theorem 9.4]{Wa} generalizes this construction to complexes over any ``exact category with duality containing $\tfrac12$.''

In the context of complexes of locally free modules over a commutative ring, an alternative local construction of a diagram similar to \eqref{eq:factor} can be found in \cite{So}, Corollaire 2.2. See also \cite{Ke}.
\end{rem}

\begin{proof}[Construction]
The cone of $\Phi$ is an acyclic complex including
\begin{equation*}
\begin{CD}
\cdots  \CC^\inv\oplus(\CC^2)^* @>f_\inv>> \CC^0\oplus (\CC^1)^* @>f_0>> (\CC^0)^*\oplus \CC^1 @>f_1>> (\CC^\inv)^*\oplus\CC^2  \cdots
\end{CD}
\end{equation*}
where
\begin{equation*}
f_\inv=f_1^*=\left(\begin{matrix}{d} &{0 } \\ {-\Phi} &{ d^*}\end{matrix}\right)
\qquad\qquad f_0=\left(\begin{matrix}{\Phi} &{ d^*} \\ {d} &{ 0}\end{matrix}\right).
\end{equation*}
Define \begin{equation}\DD^0=\(\CC^0\oplus (\CC^1)^*\)\big/\image f_\inv.
\end{equation}
 Then $f_0$ gives a symmetric isomorphism
$f_0:\DD^0\to (\DD^0)^\vee.$

Moreover, we get a chain of complexes
\begin{equation}\label{bigchain}
\begin{CD}
\CC @=[\cdots @>d>> \CC^\inv @>d>> \CC^0 @>d>> \CC^1 @>d>> \cdots]\\
  @V\alpha VV @V\Phi VV @V\Phi VV  @VbVV @| @|\\
\DD: @=[\cdots @>d^*>> (\CC^1)^* @>c>> \DD^0 @>e>> \CC^1 @>d>> \cdots]\\
 @V\tilde\Phi VV @| @| @Vf_0VV @| @| \\
\DD^* @= [\cdots @>d^*>> (\CC^1)^* @>e^*>> (\DD^0)^* @>c^*>> \CC^1 @>d>> \cdots]\\
 @V\alpha^* VV @| @| @Vb^*VV @V\Phi VV @V\Phi VV\\
\CC^* @= [\cdots @>d^*>> (\CC^1)^* @>d^*>> (\CC^0)^* @>d^*>> (\CC^\inv)^* @>d^*>> \cdots]
\end{CD}
\end{equation}
satisfying our requirements. The maps $b\colon C^0\to D^0$ and $c\colon (C^1)^*\to D^0$ are the natural ones from the definition of $\DD^0$, while $e\colon D^0\to C^1$ is induced by 
$(d,0)\colon \CC^0\oplus (\CC^1)^* \to \CC^1.$
\end{proof}

\subsection{} In our particular case, $f_{-1}=\tilde\partial$ is injective; we obtain a complex $\DD$ with $\DD^0=\bK12n/\image\tilde\partial$ and a symmetric isomorphism $f_0\colon\DD^0\to(\DD^0)^*$ inducing the form $\bq12n$. Moreover, $H^0(\DD)=H^0(\CC)=\T12n$ is the quadratic subquotient of $\DD^0=\bT12n$ by the image of $c$ in diagram \eqref{bigchain}. This completes the proof of Proposition \ref{prop:cc}.
\end{proof}

\section{Relation to Kashiwara's Construction.} \label{sec:Kash}

Let us recall the construction of Kashiwara's Maslov index, which we denote by $\tK(l_1,\ldots,l_n)$. First, for $n=3$,
 $\tK(l_1,l_2,l_3)\in\WF$ is represented by the symmetric bilinear form $\qK$ on $\TK=l_1\oplus l_2\oplus l_3$ with
$$\qK(v,w)= \frac12\big(B(v_1,w_2-w_3)+ B(v_2,w_3-w_1)+ B(v_3,w_1-w_2)\big).$$
The definition can then be extended inductively to any $n$ by the chain condition \begin{equation}\label{Kash:cocycle2}
\tK(l_1,l_2,\ldots,l_n):=\tK(l_1,l_2,\ldots,l_k)+\tK(l_1,l_k,\ldots,l_n)
\end{equation}
for any $k\in\{3,\ldots,n-1\}$. 

\begin{prop}{\label{prop:Kash} $\t12n=\tK(l_1,l_2,\ldots,l_n).$}\end{prop}

\begin{proof}{Using the chain properties \eqref{Kash:cocycle},\eqref{Kash:cocycle2}, we need only consider $n=3$.
We show that $T_{1,2,3}$ is a quadratic subquotient of $\TK$.

\begin{lemma}{$I=l_1\subset\TK$ is an isotropic subspace; $$I^\perp=\{(v_1,v_2,v_3) \mbox{ with } v_2-v_3\in l_1\}.$$}\end{lemma}

\begin{lemma}The map
 $$(v_1,v_2,v_3)\mapsto(v_2-v_3,-v_2,v_3)$$
defines an isometric surjection $I^\perp \to K_{1,2,3}$.
\end{lemma}

The lemmas, which may be checked directly, combine to give an isometric isomorphism from the non-degenerate part of $I^\perp/I$ to $T_{1,2,3}$. 
 Lemma \ref{prop:criterion} concludes the proof of Proposition \ref{prop:Kash}.}\end{proof}

\section{The Maslov Index via Sheaves.}
\label{sec:Beilinson}

\theoremstyle{remark}

In this section we verify basic properties of Beilinson's construction (see \ref{intro:Beilinson}). 

Throughout we use the following notation: if $X$ is a topological space and $W$ a vector space then $W_X$ is the constant sheaf on $X$ with fibre $W$; if $Y$ is a subset of $X$ let $j_Y$ denote its inclusion. $\mathbb D$ is the Verdier dualizing operator.

\subsection{Non-Degeneracy.} 

Since $(j_{U,!}F_U)[2]$ is the dualizing complex on $D$, the map $P\otimes P\to j_{U,!}F_U$ induced by the symplectic form on $\V$ defines a map $\phi\colon P\to(\mathbb DP)[-2]$. The fact that the bilinear form \eqref{cup} on $H^1(D,P)$ is non-degenerate follows from the following lemma. 

\begin{lemma}\label{lemma:shnondeg} The map $\phi\colon P\to(\mathbb DP)[-2]$ is an isomorphism.
\end{lemma}

\begin{proof} 

The question is local; it suffices to check that $\phi$ is an isomorphism over the open sets of the form 
$$U_i:=U\cup(i-1,i)\cup\{i\}\cup(i,i+1)\qquad\qquad i\in\ZnZ$$
since these cover $D$.
Write $M:=l_{i-1}\cap l_{i}$ and choose a decomposition
$$V= M \oplus L\oplus L'\oplus M'$$
such that $l_{i-1}= M\oplus L$ and $l_{i}=M\oplus  L'$.  Write $a,b$ for the inclusions
$$U\overset a\longrightarrow U\cup(i-1,i)\overset b\longrightarrow U_i.$$
Then on $U_i$ we have 
\begin{equation}\label{Pdecomp1}
P=b_*a_*M_U\oplus b_!a_* L_U\oplus b_*a_! L'_U\oplus b_!a_! M'_U.
\end{equation}
To compute $\mathbb DP$, note, for example, that $\mathbb Db_*a_*M_U=b_!a_!\mathbb DM_U=b_!a_!(M^*_U[2])$. Therefore, dualizing \eqref{Pdecomp1} and reversing the order of the summands, we obtain
\begin{equation}\label{Pdecomp2}(\mathbb DP)[-2]=
b_*a_* {M'}^*\oplus b_!a_*{ L'}^*_U\oplus b_*a_! L^*_U\oplus  b_!a_!M^*_U.
\end{equation}
We observe that $\phi$ is an upper-triangular matrix with respect to the decompositions \eqref{Pdecomp1}, \eqref{Pdecomp2}; moreover, the diagonal entries are isomorphisms, since the symplectic form gives isomorphisms
 $ L'\cong L^*$ and $M'\cong M^*$.
\end{proof}


\subsection{Equivalence with Algebraic Definition \ref{TheDef}.}

\def\Vert{\mathbb V\mathrm{ert}}

   The cellular cochain complex arising from the sheaf $P$ is evidently our original complex $C$ in \eqref{complex2}.  Therefore $H^1(D,P)=\T12n$. 

 As for the bilinear form, let us calculate explicitly the cup product \eqref{cup}.
We triangulate $D$ by introducing a central vertex $\infty$ and drawing radii, as in Figure 1. 

\begin{figure}[ht]
\begin{center}
\begin{picture}(0,0)%
\includegraphics{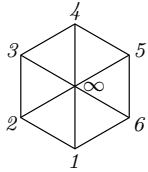}%
\end{picture}%
\setlength{\unitlength}{3947sp}%
\begingroup\makeatletter\ifx\SetFigFont\undefined%
\gdef\SetFigFont#1#2#3#4#5{%
  \reset@font\fontsize{#1}{#2pt}%
  \fontfamily{#3}\fontseries{#4}\fontshape{#5}%
  \selectfont}%
\fi\endgroup%
\begin{picture}(859,1032)(1800,-1513)
\put(2220,-1513){\makebox(0,0)[b]{\smash{{\SetFigFont{8}{9.6}{\rmdefault}{\mddefault}{\itdefault}1}}}}
\put(1862,-1273){\makebox(0,0)[rb]{\smash{{\SetFigFont{8}{9.6}{\rmdefault}{\mddefault}{\itdefault}2}}}}
\put(1869,-814){\makebox(0,0)[rb]{\smash{{\SetFigFont{8}{9.6}{\rmdefault}{\mddefault}{\itdefault}3}}}}
\put(2591,-1277){\makebox(0,0)[lb]{\smash{{\SetFigFont{8}{9.6}{\rmdefault}{\mddefault}{\itdefault}6}}}}
\put(2598,-813){\makebox(0,0)[lb]{\smash{{\SetFigFont{8}{9.6}{\rmdefault}{\mddefault}{\itdefault}5}}}}
\put(2229,-556){\makebox(0,0)[b]{\smash{{\SetFigFont{8}{9.6}{\rmdefault}{\mddefault}{\itdefault}4}}}}
\put(2280,-1026){\makebox(0,0)[lb]{\smash{{\SetFigFont{8}{9.6}{\rmdefault}{\mddefault}{\itdefault}$\infty$}}}}
\end{picture}%
\end{center}\label{fig:triangulation}
\caption{Triangulation of $D$ for $n=6$ Lagrangians.}
\end{figure}
Corresponding to this triangulation, a 1-cochain $\alpha$ with values in $P$ consists of choices of $\alpha_{\{i,i+1\}}\in l_{i}$ and $\alpha_{\{i,\infty\}}\in\V$.  
The cocycle condition amounts to
$$\alpha_{\{i+1,\infty\}}-\alpha_{\{i,\infty\}}=-\alpha_{\{i,i+1\}}.$$
A cocycle corresponds to an element $v\in\T12n$ by $v_{i}=\alpha_{\{i,i+1\}}$,  so the cocycle condition implies  $\alpha_{\{i+1,\infty\}}=-\ad{v}_{\{i,i+1\}}$ in the notation of \eqref{eq:anti}.

A 2-cocycle $\gamma$ with values in $j_{U,!}F_U$ consists of choices of $\gamma_{\{i,i+1,\infty\}}\in\V$. The evaluation $H^2(D,j_{U,!}F_U)\to F$ is given by
$\gamma\mapsto \sum_{i\in\ZnZ} \gamma_{\{i,i+1,\infty\}}.$
  Given 1-cocycles $\alpha,\beta$ the 2-cocycle $\gamma=\alpha\cup\beta$ is given by the usual formula
$$
\gamma_{\{i,i+1,\infty\}}=B(\alpha_{\{i,i+1\}}, \beta_{\{i+1,\infty\}}).
$$
All together, the pairing \eqref{cup} is given on cocycles by
$$(\alpha,\beta)\mapsto \sum_{i\in\ZnZ} B(\alpha_{\{i,i+1\}},\beta_{\{i+1,\infty\}})=\sum_{i\in\ZnZ} B(v_i,-\ad{w}_{\{i,i+1\}})=-\q12n(v,w)$$
where $v,w\in\T12n$ are the classes of $\alpha,\beta$. The last equality is the definition \eqref{eq:defq} of $\q12n$.

\begin{rem} The asymmetric formula \eqref{eq:explicit} appears if we choose an asymmetric triangulation of $D$, namely, if we draw all the diagonals from the vertex $1$.
\end{rem}

\subsection{Chain Condition.} \label{sec:BeilinsonChain}
 Let us write $(D_1,P_1)$, $(D_2,P_2)$, $(D_3,P_3)$ for the polygons and sheaves corresponding to the collections $\{l_1,\ldots,l_k\}$, $\{l_1,l_k,\ldots,l_n\}$, and $\{l_1,\ldots,l_n\}$. 

Here is a proof of the chain condition \eqref{Kash:cocycle} that is very similar to the proof of the additivity of the index of manifolds. Namely, we describe a ``bordism'' $(Y,\hat P)$ from $(D_1\sqcup D_2,P_1\oplus P_2)$ to $(D_3,P_3)$, and argue that the
image of $H^1(Y,\hat P)$ in  $H^1(D_1,P_1)\oplus H^1(D_2,P_2)\oplus H^1(\bar D_3,P_3)$ is isotropic of half the total dimension; therefore the total quadratic space is hyperbolic. Here $\bar D_3$ is $D_3$ with opposite orientation. 

In Figure 2 
 we show $Y$, with $n=4,k=3$ for concreteness.
\begin{figure}[ht]\label{fig:bordism}
\begin{center}
\begin{picture}(0,0)%
\includegraphics{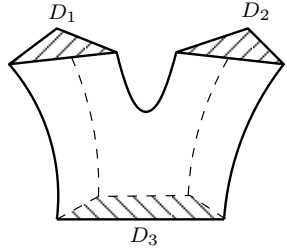}%
\end{picture}%
\setlength{\unitlength}{3947sp}%
\begingroup\makeatletter\ifx\SetFigFont\undefined%
\gdef\SetFigFont#1#2#3#4#5{%
  \reset@font\fontsize{#1}{#2pt}%
  \fontfamily{#3}\fontseries{#4}\fontshape{#5}%
  \selectfont}%
\fi\endgroup%
\begin{picture}(1786,1543)(2679,-2788)
\put(3036,-1350){\makebox(0,0)[b]{\smash{{\SetFigFont{8}{9.6}{\rmdefault}{\mddefault}{\itdefault}$D_{1}$}}}}
\put(4250,-1341){\makebox(0,0)[b]{\smash{{\SetFigFont{8}{9.6}{\rmdefault}{\mddefault}{\itdefault}$D_{2}$}}}}
\put(3551,-2752){\makebox(0,0)[b]{\smash{{\SetFigFont{8}{9.6}{\rmdefault}{\mddefault}{\itdefault}$D_{3}$}}}}
\end{picture}%
\end{center}
\caption{A ``bordism'' $Y$ between $D_1\sqcup D_2$ and $D_3$.}
\end{figure}
 The sheaf $\hat P$ on $Y$ is constructible with respect to the pictured cell decomposition, with stalks as follows: $V$ on the interior of $Y$ and on the open top and bottom faces; $l_1,l_2,l_3,l_4$ on the open front, left, back, and right faces respectively. The stalk over any other cell is the intersection of the stalks over all adjacent cells, and all the restriction maps are inclusions. 

Define $\hat D:=D_1\sqcup
D_2\sqcup\bar D_3 \subset Y$, and write $i$ for $j_{\hat D}$ and $j$ for
$j_{Y-\hat D}$. Let $\hat U$ be the interior of $Y$.
We orient $Y$ compatibly with the orientation of $\hat D$. The dualizing complex on $Y$ is then $j_{\hat U,!}F_{\hat U}[3]$. The symplectic form induces a pairing $\hat P\otimes j_{!}j_{}^!\hat P\to j_{\hat U,!}F_{\hat U}[3]$ and therefore a map $$\hat\phi\colon j_{!}j_{}^!\hat P \to (\mathbb D\hat P)[-3].$$

\begin{lemma} $\hat\phi$ is an isomorphism. In particular, it induces isomorphisms $$H^{m}(Y,j_{!}j_{}^!\hat P)\cong H^{3-m}(Y,\hat P)^*$$  for $m=0,1,2,3$.
\end{lemma}
\begin{proof} This can be proved in the same way as Lemma \ref{lemma:shnondeg}. In fact, one can reduce to Lemma \ref{lemma:shnondeg}: the pair $(Y,\hat P)$ is locally isomorphic to $(D\times I,P\boxtimes F_I)$, where $I$ is the the closed interval, and  locally $\hat\phi$ is the isomorphism
$$\hat\phi=\phi\boxtimes\mathrm{id}\colon   P\boxtimes j_{I-\partial I,!}F_{I-\partial I}\to (\mathbb DP)[-2]\boxtimes (\mathbb DF_I)[-1]=\(\mathbb D(P\boxtimes F_I)\)[-3].$$
\end{proof}

\begin{prop} The image of $\pi\colon  H^1(Y,\hat P)\to H^1(\hat D,i_{}^*\hat P)$ is isotropic of half the dimension of the total space.
\end{prop}
\begin{proof} It follows from the commutativity of the diagram
\begin{equation*}
\xymatrix{
H^1(Y,\hat P) \ar[r]^{\pi}\ar[d]_{\cong}
  & H^1(\hat D,i_{}^*\hat P) \ar[r]^{\delta}\ar[d]_{\cong}
  & H^2(Y,j_{!}j_{}^!\hat P)\ar[d]_{\cong} \\
H^2(Y,j_{!}j_{}^!\hat P)^* \ar[r]^{\delta^*} 
  & H^1(\hat D,i_{}^*\hat P)^* \ar[r]^{\pi^*} 
  & H^1(Y,\hat P)^* 
}
\end{equation*}
and the exactness of the rows.
\end{proof}

\def\Schwartz{\mathfrak S}
\def\Meas{\mathfrak M}
\def\p#1#2{p^{#1}_{#2}}
\def\dens#1#2{\left|\det {#1}\right|^{#2}}
\def\dons#1{\left|\det {#1}\right|}
\def\F{\mathcal F}
\def\Tr{\mathrm{Tr}\,}
\def\mmu#1#2{\mu_{#1,#2}^{1/2}}

\section{Relation to the Weil Representation: Introduction}
\label{sec:gamma}

From now on $F$ is a finite or local field. Fix an additive character $\psi\colon F\to \mathbb C^\times$. 
In the remainder of this article we explain how the quadratic form $\q12n$ arises in the context of the Weil representation, presenting from a new perspective the results described in \cite{LV}. 
Here is an synopsis of what follows.

\subsection{}
The Heisenberg group $G$  associated to $\V$ is $G=\V\times F$ as a set, with multiplication $$(v,s)\cdot(w,t)=(v+w,s+t+\tfrac12B(v,w)).$$ The factor $F$ is the center of $G$.
$G$ is known to have a unique irreducible unitary representation with central character $\psi$. To each Lagrangian $l_i$ is associated a realization $\rho_i$ of this representation on a Hilbert space $H(l_i)$; see \ref{sec:Hilbert}--\ref{sec:Gaction}.   There is a canonical intertwiner $\F_{j,i}\colon H(l_i)\to H(l_j)$ for any pair $l_i,l_j$; see \ref{sec:intertwiners}.
 Define $$\F_{1,2,\ldots, n}=\F_{n,n-1}\circ\cdots\circ \F_{2,1}\circ \F_{1,n}\colon H(l_n)\to H(l_n).$$
By irreducibility of $\rho_i$,  $\F_{1,2,\ldots, n}$ must be a scalar of modulus 1. On the other hand, Weil \cite{We} defined a unitary character $\gamma$ of $\WF$; see \ref{sec:Weil}. We will prove:

\begin{thm}{\label{prop:gamma} $\F_{1,2,\ldots, n}=\gamma(-\t12n)$. }\end{thm}

 This result is originally due to P. Perrin [Pe, LV], following earlier work by G. Lion \cite{Li}  and J.-M. Souriau \cite{Sou} when $F=\mathbb R$, and at the same time as similar results by R. Ranga Rao \cite{Ra}. They all consider only the case $n=3$, from which the general case follows from the  chain condition \eqref{Kash:cocycle} and the identity $\F_{i,j}=\F_{j,i}^{-1}$. In the present proof, which properly occupies section \ref{sec:ProofWeil}, we will see how our quadratic form $\q12n$ arises if we proceed directly for any $n$. 

 In fact (see Proposition \ref{preinter}) $\q12n$ logically occurs even before the introduction of the intertwiners $\F_{j,i}$.

As a corollary to the proof of Theorem \ref{prop:gamma}, we describe the self-dual measure on $\T12n$ in \ref{measure}.

\medskip

\noindent {\bf Notation.} We write $\beta$ for the projection $V\to V^*$, $\beta(x)=B(-,x).$

The unadorned symbols $\bigcap,\sum,\bigoplus,\prod,\bigotimes$ always imply the index $i\in\ZnZ$.

\subsection{The Finite Field Case}\label{fftrace} Let us summarize what happens when $F$ is a finite field. Then it makes sense to compute the trace of the scalar operator $\F(l_1,\ldots,l_n)$. 

\subsubsection{} $H(l_i)$ is the space of sections of a complex line-bundle on $V/l_i$. Representing each $\F_{i+1,i}$ by an integral kernel on $V/l_{i+1}\times V/l_i$, see \eqref{eq:F12b}, this trace takes the form
\begin{equation*}\label{thetr}\Tr\F(l_1,\ldots,l_n)=\sum\displaylimits_{\substack{v\in \bigoplus V/l_i}} 
J(v)
\end{equation*}
for some function $J$ on $\bigoplus V/l_i$.
On the other hand,  the quadratic form $q^*$ dual to $q$ is defined as in \ref{defEs} on a subspace $\Es12n\subset \bigoplus V/l_i$. 

\begin{prop}\label{p1} $J$ is supported on $\Es12n$ and
$$J(v)=|F|^{-\tfrac12\sum\dim{l_i/l_i\cap l_{i+1}}}\cdot\psi(-\tfrac12\q12n^*(v,v))\quad\forall v\in\Es12n.$$ 
\end{prop}
\begin{proof}
Compare \eqref{eq:F12b} below to \eqref{eq:dualform}. In \eqref{eq:F12b} the factor $\mmu{i}{i+1}$ can be interpreted as counting measure times $|F|^{-\tfrac12\dim l_i/l_i\cap l_{i+1}}$.
\end{proof}
Theorem  \ref{prop:gamma} now follows from the formula for $\gamma$
 recalled in \ref{ffgamma}.  

\begin{rem} The generalization of Proposition \ref{p1} to a possibly infinite field is provided by Proposition \ref{prop:traces}(ii) below.
\end{rem}

\subsubsection{} The general argument will also use what amounts to the following fact over a finite field (compare Proposition \ref{prop:traces}(i) below). 
Consider the composed function
\begin{equation*}
\begin{CD}\bigoplus l_i @>m>>
 G @>\Tr\rho_i>>\CCC$$
\end{CD}
\end{equation*}
where $m(g_1,\ldots g_n)=g_n\cdots g_2g_1$ (multiplication in $G$).

\begin{prop}\label{preinter} $\Tr\rho_i\circ m$ is supported on  $\K12n$ (see \eqref{eq:defK}) and its restriction is 
$$\Tr\rho_i\circ m(v)=|F|^{\tfrac12\dim V}\psi(\tfrac12\q12n(v,v))\qquad\forall v\in\K12n.$$
\end{prop}
\begin{proof}
First observe that if $v=(v_i)\in V^n$ then 
\begin{equation}\label{eq:mform}m(v)=\big(\sum v_i,\tfrac12\sum_{n\geq i>j\geq 1}B(v_i,v_j)\big).\end{equation}
Comparing this to \eqref{eq:explicit}, the statement now follows by computing $\Tr\rho_i$ explicitly, as in Lemma \ref{thetrace} below. 
\end{proof}

\section{Recollections about the Heisenberg Group and Weil's $\gamma$.}\label{sec:Heisgp}

A more detailed exposition of the facts contained in this section can be found in \cite{LV} \S1.2-1.4 and its appendix.

\def\OR#1#2{\Omega^\RRR_{#1}\!\(#2\)}
\def\OC#1#2{\Omega_{#1}\!\(#2\)}

\subsection{Conventions on Measures and Densities}\label{hd} 
 For $\alpha\in\RRR$, the space  of $\alpha$-densities on an $F$-vector space $X$ is the one-dimensional $\RRR$-vector space
$$\OR\alpha{X}=\{\nu\colon\det X\to\RRR\,|\,\nu(\lambda x)=|\lambda|^\alpha \nu(x),\,\forall x\in\det X,\,\lambda\in\RRR\}.$$

We identify $\OR{}{X}:=\OR{1}{X}$ with the space of real invariant measures on $X$: $\nu\in\OR{}{X}$ corresponds to the invariant measure that assigns to  $\{a_1v_1+\cdots+a_kv_k\,|\,a_i\in F, |a_i|\leq 1\}$ the volume $\nu(v_1\wedge\ldots\wedge v_k)$, for any basis $v_1,\ldots,v_k$ of $X$. 

\subsubsection{} \label{OR:aut}An isomorphism $f\colon X\to Y$ induces an isomorphism $\OR{\alpha}{X}\to\OR{\alpha}{Y}$.

 We can identify
$\OR{\alpha}{X}\otimes\OR{\beta}{X}=\OR{\alpha+\beta}{X}$
and $\OR{-\alpha}{X}=\OR{\alpha}{X}^*=\OR{\alpha}{X^*}.$
 If $Y\subset X$ then one can identify $\OR{\alpha}{X}=\OR{\alpha}{Y}\otimes\OR\alpha{X/Y}$.

 \subsubsection{} Set $\OC\alpha{X}:=\OR\alpha{X}\otimes\CCC$. Then \ref{OR:aut} works for $\Omega_{\alpha}$ as well as for $\Omega^\RRR_\alpha$.
\subsubsection{} Given $\mu\in\OR{}{X}$ there is a unique `dual' measure $\mu^*\in\OR{}{X^*}$ with $\<\mu,\mu^*\>=1$. Given an isomorphism $s\colon X\to X^*$, there is a unique positive measure $\mu_X\in\OR{}{X}$ that is `self-dual,' i.e. $s_*(\mu_X)=\mu_X^*$.

\subsubsection{Half-Densities}
    For $\nu\in\OR{}{X}$ a positive measure, define $\nu^{1/2}\in\OR{1/2}{X}$ by $\nu^{1/2}(x):=\left|\nu(x)\right|^{1/2}.$ Then $\nu^{1/2}\otimes\nu^{1/2}=\nu$ using $\OR{1/2}{X}\otimes\OR{1/2}{X}=\OR{}{X}$.

\subsection{The Representation Spaces $H(l_i)$}\label{sec:Hilbert}
Let $H_s(l_i)$ be the space of functions
 $\phi\colon\V\to \OC{1/2}{V/l_i}$  satisfying
\begin{equation}\label{transform}
\phi(x+a)= \[x,a\]\cdot\phi(x) \quad \mbox{ $\forall x\in V,a\in l_i$}
\end{equation}
and such that $\phi$ is Schwartz modulo $l_i$ (i.e. $\phi$ is smooth, and $|\phi(x)|,$ which is constant along $l_i$, decays rapidly as a function on $V/l_i$).

Let $H(l_i)$ be the completion of $H_s(l_i)$ 
 with respect to the following norm:  $\phi\bar\phi$ is an $l_i$-invariant function on $\V$ with values in $\OC{}{V/l_i}$; the norm of $\phi$ is defined by
$$|\phi|^2=\int_{V/l_i} \phi\bar\phi.$$

\def\pipi{{\p{i'}i}}
\def\piip{\p{i}{i'}}
\def\pnpn{{\p{n'}n}}
\def\pnnp{\p{n}{n'}}

\subsection{The Representation of $G$ on $H(l_i)$}\label{sec:Gaction}
$H(l_i)$ is the space of the representation $\rho_i$ of $G$ induced from the character $(a,t)\mapsto\psi(t)$ of $l_i\oplus F\subset G$. Explicitly, for $\phi\in H(l_i)$ and $(v,t)\in G$, 
\begin{equation}\label{eq:Gaction}
\rho_i(v,t)\phi(x)=\phi(x-v)\cdot\[v,x\]\cdot\psi(t).
\end{equation}


\subsubsection{The Character}\label{sec:haction}

Suppose $h$ is a $C_0^\infty$ measure on $G$ (i.e. $h$ is a smooth, compactly supported function times a Haar measure). It gives rise to an operator $\rho_i(h)$ on $H(l_i)$ via 
\begin{equation*}
\rho_i(h)\phi(x):=\int\displaylimits_{(v,t)\in G} \rho_i(v,t)\phi(x) \cdot h(v,t)
\end{equation*}
or, explicitly,
\begin{equation}\label{eq:haction}
\rho_i(h)\phi(x)=\int\displaylimits_{(v,t)\in G}\phi(x-v)\cdot\[v,x\]\cdot\psi(t)\cdot h(v,t).
\end{equation}

The operator $\rho_i(h)$ is always trace-class (see the proof of Lemma \ref{thetrace} below). One defines the character 
$\Tr\rho_i$ of $\rho_i$ 
as a generalized function on $G$, i.e. as a functional on the space of $C_0^\infty$  measures, by 
\begin{equation}\label{trdef}\<\Tr\rho_i,h\>:=\Tr \rho_i(h).
\end{equation}

\begin{lemma}\label{thetrace}  $\Tr\rho_i=\delta\otimes\psi$ as a generalized function on $V\times F$. 
Here $\delta$ is the delta-function at $0\in V$  (i.e. the delta-measure divided by self-dual Haar measure).  
\end{lemma}

\begin{proof}
Changing variables $v\mapsto x-v$ in \eqref{eq:haction},
\begin{equation*}
\begin{aligned}\rho_i(h)\phi(x)&=
\int\displaylimits_{(v,t)\in G} 
\phi(v)\cdot\[x,v\]\cdot\psi(t)\cdot h(x-v,t)\\
&=\int\displaylimits_{b\in V/l_i}
\phi(b)
\int\displaylimits_{(a,t)\in l_i\oplus F} \[x+ b,a+b\]\cdot\psi(t)\cdot h(x-a-b,t).
\end{aligned}
\end{equation*}
For the last formula we have used \eqref{transform}.
The inner integral is smooth in $(x,b)\in V/l_i\times V/l_i,$ so we obtain
\begin{equation*}
\Tr\rho_i(h)=\int\displaylimits_{b\in V/l_i}\int\displaylimits_{(a,t)\in l_i\oplus F} \[2b,a\]\cdot\psi(t)\cdot h(-a,t).
\end{equation*}
The integrals over $a$ and $b$ are Fourier inverse, leaving exactly 
\begin{equation*}\label{eq:character}
\Tr\rho_i(h)=\int\displaylimits_{t\in  F}  h(0,t)\cdot\psi(t)\cdot\mu_V^*.
\end{equation*}
\end{proof}

\subsection{Intertwiners}\label{sec:intertwiners}
For each pair of Lagrangians $l_i,l_j$ there is a canonical unitary intertwiner
$\F_{j,i}\colon H(l_i)\to H(l_j)$.  According to \cite{LV} it is determined by its effect on $\phi\in H_s(l_i)$, which is
\begin{equation}\label{eq:F12}
\F_{j,i}\phi(y):=\int\displaylimits_{x \in l_j/(l_i\cap l_j)} \phi(x+y)\cdot\[x,y\]\cdot \mmu{i}{j}\quad \in H_s(l_j).
\end{equation}
Here $\mu_{i,j}\in\OR{}{(l_i+l_j)/(l_i\cap l_j)}\cong\OR{}{l_i}\otimes\OR{}{l_j}\otimes(\OR{}{l_i\cap l_j}^*)^{\otimes 2}$ is the self-dual measure on the symplectic space  $(l_i+l_j)/(l_i\cap l_j)$, so that we have
\begin{equation*}\begin{aligned}
\phi(x+y)\cdot\mu_{i,j}^{1/2}\in &\OC{1/2}{V/l_i}\otimes\OC{1/2}{l_i}\otimes\OC{1/2}{l_j}\otimes\OC{}{l_i\cap l_j}^*\\
&\cong\OC{1/2}{V/l_j}\otimes\OC{}{l_j/(l_i\cap l_j)}.
\end{aligned}
\end{equation*}

The next lemma represents $\F_{{i+1},i}$ by an integral kernel.

\begin{lemma} For $\phi\in H_s(l_i)$ and any $y\in V$
\begin{equation}\label{eq:F12b}
\F_{{i+1},i}\phi(y)=\!\!\!\!\!\int
\displaylimits_{\substack{x\in V/l_i:\\x-y\in l_i+l_{i+1}}}
\!\!\!\!\!\!\!\!\!\! \phi(x)\cdot\psi\(\tfrac12\<\rii\(x,y\),x-y\>\)\cdot
\mmu{i}{i+1}.
\end{equation}
(For the definition of $\rii$, see \ref{sec:r}.)
\end{lemma}

\begin{proof} Choose a complement $A$ to $l_i\cap l_{i+1}$ in $l_{i+1}$; make the change of variable $z=x+y$ in the definition \eqref{eq:F12} to get
\begin{equation*}
\F_{i+1,i}\phi(y)=\int\displaylimits_{\substack{z\in V:\\z-y\in A}}\phi(z)\cdot\[z,y\]\cdot \mmu{i}{i+1}.
\end{equation*}
We observe that if $z-y\in A\subset l_{i+1}$ then $B(z,y)=\<\rii(z,y),y-z\>$.  With this substitution, the integrand only depends on $z$ modulo $l_i$. Finally, 
projection along $l_i$ gives a bijection $\{z\in V: z-y\in A\}\to\{x\in V/l_i: x\equiv y\bmod l_i+l_{i+1}\}$.
\end{proof}

\subsection{Weil's Character $\gamma$}\label{sec:gammadef}
\label{sec:Weil}
Suppose $(T,q)$ is a quadratic space.  Let $f\mapsto  f^\wedge$ denote the Fourier transform with respect to $\psi$.  Consider the functions
$$f_q\colon x\mapsto\psi\(\tfrac12q(x,x)\) \qquad\qquad
f_{-q^*}:x\mapsto\psi\(-\tfrac12q^*(x,x)\)$$
on $T$ and $T^*$ respectively. 
 Weil \cite{We} shows that there exists a character $\gamma$ of $W(F)$ such that the measure $f_q^\wedge$ on $T^\vee$ is represented by 
\begin{equation}\label{gammadef}
f_q^\wedge=\gamma(q)\cdot f_{-q^*}\,\mu_{T}^*
\end{equation}
where $\mu^*_{T}$ is the self-dual measure on $T^*$.

\subsubsection{Finite Field Case} \label{ffgamma}
Suppose that $F$ is a finite field, and also allow $q$ to degenerate as a quadratic form on $T$. Then one finds
$$\gamma(q)=|F|^{-\tfrac12(\dim T+\dim\ker q)}\sum_{x\in T}\psi(-\tfrac12q(x,x)).$$

\section{Proof of Theorem \ref{prop:gamma}.}\label{sec:ProofWeil}\label{sec:final}

\subsection{} Generalizing the argument of \ref{fftrace}, we compute the trace of the operator
$$\sigma(g_1,\ldots,g_n):= \rho_n(g_n)\circ\F_{n,n-1}\circ\cdots\circ\rho_1(g_1)\circ\F_{1,n}\colon H(l_n)\to H(l_n)$$
as a generalized function in $(g_i)\in\bigoplus l_i$. Since $\F_{i+1,i}$ intertwines $\rho_i$ and $\rho_{i+1}$, 
$$\sigma(g_1,\ldots,g_n)= \F_{1,2,\ldots,n}\cdot\rho_n(g_n\cdots g_1).$$
We obtain $\F_{1,2,\ldots,n}$ by dividing $\Tr\sigma$ by the trace of 
$\rho_n\circ m$, where $m(g_1,\ldots,g_n)=g_n\cdots g_1\in G$. 

  These traces make sense under the conditions of the following lemma. 

\begin{lemma}\label{makesense} The image $m(\bigoplus l_i)\subset G$ lies transverse to the support of $\Tr\rho_n$ if \begin{equation}\label{assumption}
\sum l_i=V\qquad\qquad\mbox{(equivalently, $\bigcap l_i=0$)}. 
\end{equation}
\end{lemma}

\begin{proof}[Proof of Lemma \ref{makesense}] 
 Assuming \eqref{assumption},
 $m$ is a submersion onto $G$ even when restricted to $\bigoplus (l_i\oplus F)$. 
 It remains to note from Lemma \ref{thetrace} that $\Tr\rho_n$ is equivariant under the action of $m(F^n)=F\subset G$. 
\end{proof}

\subsection{} In proving Theorem \ref{prop:gamma}
we may always assume condition \eqref{assumption}. Indeed, one sees from the definition \eqref{eq:F12} of the intertwiners that the number $\F_{1,2,\ldots,n}$ does not change if we replace $V$ by $(\sum l_i)/\bigcap l_i$ and each $l_j$ by $l_j/\bigcap l_i$; nor is our quadratic space affected.

\subsection{} Now we describe two generalized functions $Q$ and $Q'$ on $\bigoplus l_i$, which we will  identify with $\Tr\rho_n\circ m$ and $\Tr\sigma$ in Proposition \ref{prop:traces}.  We use the following objects: 

\medskip
\begin{tabular}{ll}
$f_q\colon x\mapsto\psi\(\tfrac12\q12n(x,x)\)$ & Function on $\T12n$\\
$f_{-q^*}\colon x\mapsto\psi\(-\tfrac12\q12n^*(x,x)\)$ & Function on $\T12n^*$\\
$\mu_{T},\mu_V$ & Self-dual measures on $\T12n$, $V$ \\
$1$ & Constant function on $\bigoplus l_i\cap l_{i+1}$ \\
$\delta=(\mu_{V}^*)^\wedge$ & Delta function at $0\in V$ \\
$\delta'=1^\wedge$ & Delta measure at $0\in \bigoplus (l_i\cap l_{i+1})^*$ 
\end{tabular}
\smallskip

\noindent
Referring to the complex $C$ \eqref{complex2}, we see that $1\otimes f_q$ is a function on $\ker\Sigma$ and that therefore, assuming \eqref{assumption},
$$Q:=1\otimes f_q  \otimes \delta$$ 
defines a generalized function on $\bigoplus l_i$.
Similarly, 
$$Q'^\wedge:= \delta'\otimes f_{-q^*}\,\mu_{T}^*\otimes\mu_{V}^*$$
defines a measure on $\bigoplus l_i^*$; let $Q'$ be its Fourier transform. 
  From \eqref{gammadef}  we find 
\begin{equation}\label{newgammadef}
Q'=\gamma(-\t12n)\cdot Q.\end{equation}
Thus Theorem \ref{prop:gamma} follows from:

\begin{prop} \label{prop:traces} Assume $\eqref{assumption}$. Then, as generalized functions on $\bigoplus l_i$, 
\begin{itemize}
\item[(i)] $\Tr\rho_n\circ m=Q$.
\item[(ii)] $\Tr\sigma=Q'$.
\end{itemize}
\end{prop}
\begin{proof}
 The first statement follows from the calculation of $\Tr\rho_i$ in Lemma \ref{thetrace} and comparing \eqref{eq:mform} to the formula \eqref{eq:explicit} for $\q12n$. Explicitly, 
\begin{equation}\label{finaltrrho}
\begin{aligned}
\<\Tr\rho,h\> = \int\displaylimits_{\substack{v\in\bigoplus l_i \\ \sum v_i=0}}
\psi\big(\tfrac12\sum_{n\geq i>j\geq1} B(v_i,v_j)\big) \cdot h(v)\cdot\mu_V^*.
\end{aligned}
\end{equation}

For the second statement, note that, restricting $\rho_i$ to measures $h$ supported on $l_i\subset G$, and using \eqref{transform}, 
the formula \eqref{eq:haction} for $\rho_i(h)$ becomes simply
$$\rho_i(h)\phi(x)=\phi(x)\cdot  h^\wedge(\beta(x)).$$
Therefore, using $\eqref{eq:F12b}$, we find
\begin{equation}\label{finaltrphi}
\begin{aligned}
\<\Tr\sigma,h\>&=
\!\!\!\!\!\!\!\!\!\!\!\!\!\!\!\!\!\!\displaystyle\int\displaylimits_{\substack{x\in\bigoplus V/l_i \\ x_{i}\equiv x_{i+1}\bmod l_i+l_{i+1}}}\!\!\!\!\!\!\!\!\!\!\!\!\!\!\!\!\!
\prod\big(\psi(\tfrac12 \<\rii(x_i,x_{i+1}),x_{i}-x_{i+1}\>)\cdot\mmu{i}{i+1}\big)\cdot
 h^\wedge(\beta(x))
\end{aligned}
\end{equation}
which, compared to \eqref{eq:dualform}, is $\<Q'^\wedge,  h^\wedge\>=\<Q', h\>$ up to positive scale.

Since $\F_{1,2,\ldots,n}$ and $\gamma(-\t12n)$ both have modulus 1, we may conclude that 
$\<\Tr\sigma,h\>=\<Q', h\>$ on the nose.
\end{proof}
\subsection{The self-dual measure on $\T12n$}\label{measure}
The expressions \eqref{finaltrrho}, \eqref{finaltrphi} for $Q$ and $Q'$ allow us to identify explicitly the self-dual measure $\mu_{T}$ on $\T12n$. It can be described as follows. 

Using the complex $C$ from \eqref{complex2}, 
we can identify
\begin{equation*}\OC{}{\T12n}=
\OC{}{\bigoplus l_i\cap l_{i+1}}^*\otimes\OC{}{\bigoplus l_i} \otimes
\OC{}{V}^*\otimes\OC{}{ \bigcap l_i\oplus V\big/\sum l_i)}
.\end{equation*}
Now, $\OC{}{V}^*$ has a canonical element $\mu_V^*$ corresponding to self-dual measure on $V$, and similarly for $\OC{}{\bigcap l_i\oplus V\big/\sum l_i)}$. 
On the other hand
$$\OC{}{\bigoplus l_i\cap l_{i+1}}^*\otimes\OC{}{\bigoplus l_i}=
\bigotimes\OC{1/2}{(l_i+l_{i+1})/(l_i\cap l_{i+1})}$$
each factor of which again has a self-dual element $\mu_{i,i+1}^{1/2}$.

The product of all these canonical elements is $\mu_T$.


\begin{thebibliography}{CLM}
\bibitem[CLM]{CLM} S. Cappell, R. Lee, and E. Miller. ``On the Maslov Index.'' Comm. Pure Appl. Math. {\bf 47} (1994), 121--186.
\bibitem[KS]{KS} M. Kashiwara and P. Schapira. {\it Sheaves on Manifolds.} Berlin; New York: Springer, 1990.
\bibitem[Ke]{Ke} G. Kempf. ``Deformations of Semi-Euler Characteristics.'' American J. of Math. {\bf 114} (1992), 973--978.
\bibitem[Li]{Li} G. Lion. ``Indice de Maslov et repr\'esentation de Weil.'' {\it Trois textes sur les representations des groupes nilpotents et resolubles.} Publications Math\'ema\-tiques de l'Universit\'e Paris VII (1978).
\bibitem[LV]{LV} G. Lion and M. Vergne, {\it The Weil representation, Maslov index and Theta series.} Progress in Math 6. Boston: Birkha\" user, 1980.
\bibitem[Pe]{Pe} P. Perrin. ``Repr\'esentations de Schr\"odinger, indice de Maslov et groupe metaplectique.'' {\it Non Commutative Harmonic Analysis and Lie Groups (Marseille-Luminy, 1980),} 370--407. Lecture Notes in Math 880. Berlin; New York: Springer, 1981.
\bibitem[Ra]{Ra} R. Ranga Rao. ``On some explicit formulas in the theory of Weil representations.'' Pacific J. Math. {\bf 157} (1993), 335--371.
\bibitem[So]{So} C. Sorger. ``La semi-caract\'eristique d'Euler-Poincar\'e des faisceaux $\omega$-quad\-ratiques sur un sch\'ema de Cohen-Macaulay.'' Bull. Soc. Math. France {\bf 122} (1994), 225--233.
\bibitem[Sou]{Sou} J.-M. Souriau. ``Construction explicite de l'indice de Maslov, applications.'' {\it  Group Theoretical Methods in Physics (Fourth International Colloquium, Nijmegen, 1975),} 117-148. Lecture Notes in Physics 50. Berlin; New York: Springer, 1976.
\bibitem[Th]{Th} T. Thomas. ``The character of the Weil representation.'' http://www.arxiv.org/math/ 0610644/.
\bibitem[Wa]{Wa} C. Walter. ``Gro\-then\-dieck-Witt groups of tri\-an\-gulated cate\-gories.'' Pre\-print, July 1, 2003, K-theory Pre\-print Ar\-chives, http://www.math.uiuc.edu/K-theory/0589/. 
\bibitem[We]{We} A. Weil. ``Sur certains groupes d'op\'erateurs unitaires.'' {\it Acta Math.} {\bf 111} (1976), 143--211.
\end{thebibliography}
\end{document}